\documentclass[11pt]{article}
\usepackage[utf8]{inputenc}
\usepackage[T1]{fontenc}
\usepackage{lmodern}
\usepackage{amsmath,amssymb,amsthm,mathtools}
\usepackage{enumitem}
\usepackage{geometry}
\usepackage{xcolor}
\usepackage{microtype}
\usepackage{hyperref}
\hypersetup{colorlinks=true,linkcolor=blue!55!black,
 citecolor=blue!55!black,urlcolor=blue!55!black,
 pdftitle={Recurrence and Transience in Monotone Environments: Increasing Conductances and Growing Domains},
 pdfauthor={Nordine Anis Moumeni}}
\allowdisplaybreaks[1]
\newtheorem{theorem}{Theorem}[section]
\newtheorem{proposition}[theorem]{Proposition}
\newtheorem{lemma}[theorem]{Lemma}
\newtheorem{corollary}[theorem]{Corollary}

\newtheorem{remark}[theorem]{Remark}

\newcommand{\N}{\mathbb N}
\newcommand{\Z}{\mathbb Z}
\newcommand{\R}{\mathbb R}
\newcommand{\PP}{\mathbb P}
\newcommand{\EE}{\mathbb E}
\newcommand{\1}{\mathbf 1}
\newcommand{\cF}{\mathcal F}
\newcommand{\cE}{\mathcal E}
\newcommand{\eps}{\varepsilon}
\newcommand{\Reff}{R_{\mathrm{eff}}}
\newcommand{\supp}{\operatorname{supp}}
\newcommand{\Capacity}{\operatorname{Cap}}

\title{Recurrence and Transience in Monotone Environments\\
Increasing Conductances and Growing Domains}
\author{Nordine Anis Moumeni}
\date{November 2026}

\begin{document}
\maketitle
\begin{abstract}
We study discrete-time random walks on a locally finite graph whose
conductances are non-decreasing in time and bounded on each edge.
We compare its hitting probabilities with those of the homogeneous walk
whose conductances are fixed at the starting time. These estimates give
another proof of two results of Li and Park, who work under more general
assumptions. If the upper network is recurrent, the walk visits every
vertex infinitely often almost surely. If the initial network is
transient, the expected number of visits to each vertex is finite.

Let us first consider the probability of leaving a finite set before
hitting a given vertex. The difference between this probability and
its value for the frozen walk is bounded in absolute value by the
square root of the product of a capacity increment and a dual energy.
On an initially transient network, we obtain a lower bound on the
probability of never hitting the vertex. Both proofs use a comparison
of two Dirichlet forms and a contraction property of the killed
transition operators.

We also consider walks on increasing finite sets. On each prescribed
time interval, called a plateau, the transition operator is fixed.
A lower bound on transition probabilities after mixing and an upper
bound on the expected number of returns give a sufficient condition
for every vertex to be visited infinitely often almost surely.
For uniformly elliptic lazy walks on boxes of radius $R_n$ in $\Z^d$,
$d\geq3$, it is enough that $\sum_n(L_n-b_n)_+/R_n^d$ diverges, where
$L_n$ is the duration of the $n$th plateau and
$b_n=\lceil aR_n^2\rceil$ for a constant $a$ depending only on the
dimension and the ellipticity bound.
\end{abstract}
\newpage
\tableofcontents
\newpage
\section{Introduction and main results}\label{sec:introduction}

Let $G=(V,E)$ be a connected, locally finite graph with a countably infinite
vertex set. We consider a deterministic sequence $(c_t)_{t\geq0}$ of positive
symmetric conductances on $E$. We set $c_t(x,y)=0$ for nonadjacent
vertices. For every $t\geq0$ and every $x,y\in V$, we write
\[
 \pi_t(x)=\sum_{z:z\sim x}c_t(x,z)
 \qquad\text{and}\qquad
 K_t(x,y)=\frac{c_t(x,y)}{\pi_t(x)}\1_{\{x\sim y\}}.
\]
For $s\geq0$ and $x\in V$, let $(X_t)_{t\geq s}$ be the Markov chain
started from $x$ at time $s$, with transition
operator $K_t$ from time $t$ to time $t+1$. We denote its law and expectation by $\PP_{s,x}$ and
$\EE_{s,x}$. For $s\leq t$, we write
\[
 K_{s,t}=K_s\cdots K_{t-1}
 \qquad\text{with}\qquad K_{s,s}=I.
\]
For a subset $A\subseteq V$, we denote by $\tau_A^s$ the first hitting time
of $A$ at or after time $s$:
\[
 \tau_A^s=\inf\{t\geq s:X_t\in A\}.
\]
For a vertex $o\in V$, we write $\tau_o^s=\tau_{\{o\}}^s$.

We assume that the conductances are non-decreasing and bounded on each
edge. More precisely, there is a conductance configuration $\bar c$ such that
\begin{equation}\label{eq:monotone}
 \forall t\geq0,\ \forall e\in E,\qquad
 0<c_0(e)\leq c_t(e)\leq c_{t+1}(e)\leq\bar c(e)<\infty.
\end{equation}
Thus, for every $e\in E$, the limit $c_\infty(e)=\lim_{t\to\infty}c_t(e)$
exists and is finite. Since the graph is locally finite, one also has
\[
 \forall x\in V,\qquad
 \pi_\infty(x):=\sum_{y:y\sim x}c_\infty(x,y)<\infty.
\]

In the time-homogeneous case, recurrence is characterized by effective
resistance \cite[Chapter~2]{LP16book}. If $(D_n)_{n\geq1}$ is an
increasing exhaustion of $V$ by finite connected sets containing $o$, then
\[
 (G,c)\text{ is recurrent}
 \quad\Longleftrightarrow\quad
 \lim_{n\to\infty}\Reff^c(o,D_n^c)=\infty.
\]
In particular, if $c\leq\bar c$ and $(G,\bar c)$ is recurrent, then
$(G,c)$ is recurrent. Indeed, Rayleigh monotonicity gives
\[
 \forall n\geq1,\qquad
 \Reff^c(o,D_n^c)\geq\Reff^{\bar c}(o,D_n^c).
\]
Thus, under this assumption, the walk obtained by fixing the conductances
at any one time is recurrent. The question is whether the
time-inhomogeneous walk is also recurrent when the transition operator
changes at every step.

In the first part of this paper, we compare the hitting probabilities of
these two walks, fixing the conductances at the starting time. We first
bound the absolute difference of the exit probabilities in a finite
domain. On an infinite graph with a transient initial network, we obtain
a lower bound on the probability of never hitting a given vertex.
From these estimates we recover two results of Li and Park~\cite{LP26}.
If the upper network is recurrent, the walk visits every vertex infinitely
often almost surely. If the initial network is transient, the expected
number of visits to each vertex is finite. We begin with the recurrence
statement.

\begin{theorem}[Recurrence]\label{thm:recurrence}
Let $(c_t)_{t\geq0}$ be a deterministic sequence of conductances
satisfying~\eqref{eq:monotone}. Assume that $(G,\bar c)$ is recurrent.
Then, for every $s\geq0$ and every $x,o\in V$, we have
\[
 \PP_{s,x}(\tau_o^s<\infty)=1.
\]
Moreover, from every starting point and time, every vertex is visited
infinitely often almost surely.
\end{theorem}

This is the conclusion of \cite[Theorem~1.6]{LP26} under our assumptions.
The bounds $c_0(e)$ and $\bar c(e)$ may depend on $e$. In particular,
we do not assume uniform ellipticity. To prove
Theorem~\ref{thm:recurrence}, we will use the following comparison estimate.

Let $c$ be a conductance configuration. For functions of finite energy,
we denote its Dirichlet form by
\[
 \cE_c(f,g)=\frac12\sum_{x,y\in V}c(x,y)
 (f(x)-f(y))(g(x)-g(y)),
\]
where $c(x,y)=0$ when $x$ and $y$ are not adjacent. We write
$\cE_t=\cE_{c_t}$. Let $D$ be a finite connected proper subset of $V$
containing $o$, and set
\[
 F=D\setminus\{o\}\qquad\text{and}\qquad A=\{o\}\cup D^c.
\]
For every $t\geq0$, let $u_t^D$ be the function which satisfies
\[
 u_t^D(o)=0,\qquad u_t^D=1\text{ on }D^c,
 \qquad \cE_t(u_t^D,f)=0\quad\text{for every }f|_A=0.
\]
Thus $u_t^D(x)$ is the probability of exiting $D$ before hitting $o$
for the time-homogeneous walk with conductances $c_t$.
For every $t\geq0$, we define
\[
 C_t^D=\cE_t(u_t^D,u_t^D)=\frac1{\Reff^{c_t}(o,D^c)}.
\]
By Rayleigh monotonicity, $(C_t^D)_{t\geq0}$ is non-decreasing and bounded.
We write
\[
 C_*^D=\lim_{t\to\infty}C_t^D
 \qquad\text{and}\qquad
 C_{\bar c}^D=\frac1{\Reff^{\bar c}(o,D^c)}.
\]
In particular, $C_t^D\leq C_*^D\leq C_{\bar c}^D$ for every $t\geq0$.
For a signed measure $\nu$ on $F$, define
\[
 H_t^D(\nu)=\sup_{f|_A=0}\{2\nu(f)-\cE_t(f,f)\}.
\]
This supremum is finite since the Dirichlet form is positive definite
on the space of functions vanishing on $A$. For $x\in F$, we write
\[
 R_t^D(x,A):=H_t^D(\delta_x)
\]
for the effective resistance between $x$ and the set $A$. Here the
vertices of $A$ are identified and held at potential zero.

Let $K_t^F$ be the restriction of $K_t$ to $F\times F$. It is the
transition operator of the walk killed upon entering $A$. For $s\geq0$
and $x\in F$, we write
\[
 p_s^D(x)=\PP_{s,x}(\tau_{D^c}^s<\tau_o^s).
\]
We can now compare the two walks started from $x$ at time $s$.
The first bound keeps track of the distribution after one step of
the killed walk.

\begin{theorem}[Comparison with the frozen network]
\label{thm:finite-stability}
Let $(c_t)_{t\geq0}$ satisfy~\eqref{eq:monotone}. Let $o\in V$ and let
$D\subset V$ be a finite connected proper subset containing $o$.
Then, for every $s\geq0$ and every $x\in D\setminus\{o\}$, we have
\begin{align}
 |p_s^D(x)-u_s^D(x)|
 &\leq\sqrt{(C_*^D-C_s^D)\,
                 H_s^D(\delta_xK_s^F)}\label{eq:finite-stability}\\
 &\leq\sqrt{(C_{\bar c}^D-C_s^D)\,R_s^D(x,A)}.
 \nonumber
\end{align}
\end{theorem}

In particular, if $c_t=c_s$ for every $t\geq s$, then $C_*^D=C_s^D$
and the estimate gives $p_s^D=u_s^D$. In general, the error is controlled
by the capacity difference $C_*^D-C_s^D$. The increments
$c_{t+1}(e)-c_t(e)$ need not satisfy any further bound
beyond~\eqref{eq:monotone}.

We will also use the
free resistance between $x$ and $o$, defined by
\begin{equation}\label{eq:free-resistance}
 R_t^{\mathrm{free}}(x,o)=
 \sup_{0<\cE_t(f,f)<\infty}
 \frac{|f(x)-f(o)|^2}{\cE_t(f,f)},
\end{equation}
where the supremum is over all real functions of finite energy.
This resistance is finite. Indeed, summing the differences of $f$ along
a simple path from $x$ to $o$ and applying Cauchy--Schwarz bounds
$R_t^{\mathrm{free}}(x,o)$ by the sum of the edge resistances along that path.
Moreover, one has $R_t^D(x,A)\leq R_t^{\mathrm{free}}(x,o)$.

\begin{corollary}[Escape estimate]\label{cor:escape}
Let $(c_t)_{t\geq0}$ satisfy~\eqref{eq:monotone}. Let $o\in V$ and let
$D\subset V$ be a finite connected proper subset containing $o$.
Then, for every $s\geq0$ and every $x\in D\setminus\{o\}$, we have
\begin{equation}\label{eq:escape}
 \PP_{s,x}(\tau_{D^c}^s<\tau_o^s)
 \leq 2\sqrt{
 \frac{R_s^{\mathrm{free}}(x,o)}{\Reff^{\bar c}(o,D^c)}}.
\end{equation}
\end{corollary}

Suppose that $(G,\bar c)$ is recurrent and let $(D_n)_{n\geq1}$ be an
increasing exhaustion by finite connected sets containing $x$ and $o$.
The numerator in~\eqref{eq:escape} is finite and does not depend on $n$,
whereas $\Reff^{\bar c}(o,D_n^c)\to\infty$. Hence
\[
 \lim_{n\to\infty}\PP_{s,x}(\tau_{D_n^c}^s<\tau_o^s)=0.
\]
We will prove that the walk almost surely exits each $D_n$ or hits $o$
in finite time. The preceding limit then implies
$\PP_{s,x}(\tau_o^s<\infty)=1$, as stated in
Theorem~\ref{thm:recurrence}.

\subsection{Transience}

We now assume that the initial network $(G,c_0)$ is transient.
By Rayleigh monotonicity, $(G,c_t)$ is transient for every $t\geq0$.
Let $o\in V$. We denote by $C_c(V)$ the functions on $V$ with finite
support. For every $t\geq0$, define
\[
 \kappa_t(o)=\Capacity_{c_t}(\{o\})
 =\inf\{\cE_t(f,f):f\in C_c(V),\ f(o)=1\}.
\]
The sequence $(\kappa_t(o))_{t\geq0}$ is non-decreasing and satisfies
\[
 \forall t\geq0,\qquad
 0<\kappa_0(o)\leq\kappa_t(o)\leq\pi_t(o)\leq\pi_\infty(o)<\infty.
\]
Thus its limit exists, and we write
$\kappa_*(o)=\lim_{t\to\infty}\kappa_t(o)$.
Let $\PP_x^{c_t}$ be the law of the time-homogeneous walk with
conductances $c_t$, started from $x$. We set
\[
 u_t(x)=\PP_x^{c_t}(\tau_o=\infty).
\]
For every finitely supported signed measure $\nu$, we also define
\[
 H_t^o(\nu)=
 \sup_{\substack{f\in C_c(V)\\ f(o)=0}}
 \{2\nu(f)-\cE_t(f,f)\}.
\]
Recall that~\eqref{eq:free-resistance} allows all functions of finite
energy. Here the test functions vanish at $o$ and have finite support.
This is the condition we obtain by grounding the exterior of a finite
domain and then taking an exhaustion.

\begin{theorem}[Transience and late-time escape]\label{thm:transience}
Let $(c_t)_{t\geq0}$ satisfy~\eqref{eq:monotone}. Assume that $(G,c_0)$
is transient and let $o\in V$. Then, for every $s\geq0$ and every
$x\in V\setminus\{o\}$, we have
\begin{equation}\label{eq:transience-stability}
 \PP_{s,x}(\tau_o^s=\infty)
 \geq u_s(x)-
 \sqrt{(\kappa_*(o)-\kappa_s(o))H_s^o(\delta_x)}.
\end{equation}
For a walk started at $o$ at time $t$, let
\[
 q_t(o)=\PP_{t,o}(X_n\ne o\text{ for every }n>t).
\]
Then the probability of never returning to $o$ satisfies
\begin{equation}\label{eq:last-exit}
 \liminf_{t\to\infty}q_t(o)
 \geq\frac{\kappa_*(o)}{\pi_\infty(o)}>0.
\end{equation}
Consequently, from every deterministic starting point and time,
every vertex is visited finitely often almost surely, and the expected
number of visits to each vertex is finite.
\end{theorem}

Li and Park~\cite[Theorem~1.8]{LP26} proved transience and gave a
quantitative bound on the expected number of visits. Our comparison
with $u_s$ in~\eqref{eq:transience-stability} depends on the remaining
capacity increment $\kappa_*(o)-\kappa_s(o)$. We only use the existence
of $\kappa_*(o)$, without assuming that
$\kappa_*(o)=\Capacity_{c_\infty}(\{o\})$.

\begin{remark}[Nonadaptive random environments]\label{rem:nonadaptive}
The recurrence and transience results also apply to random conductances
under the following assumption. Let
$\mathcal G_\infty=\sigma(c_t:t\geq0)$ and suppose that
\begin{equation}\label{eq:quenched}
 \PP(X_{t+1}=y\mid\mathcal G_\infty\vee\sigma(X_0,\ldots,X_t))
 =K_t(X_t,y)\quad\text{almost surely}.
\end{equation}
Assume that the corresponding hypotheses on the conductances hold
almost surely. Conditionally on $\mathcal G_\infty$, the environment is
fixed and the walk has transition operators $(K_t)_{t\geq0}$. We can
therefore apply the deterministic theorems. Under the transience
assumptions, the conditional expected number of visits to each vertex is
finite for almost every environment. Finiteness of the expected number of visits after
averaging over environments requires integrability of this conditional
expectation.
\end{remark}

\subsection{Recurrence on growing finite sets}

We now consider a walk on an increasing sequence of finite sets.
Let $V_1\subseteq V_2\subseteq\cdots$ be finite sets, and write
$V_\infty=\bigcup_{n\geq1}V_n$. For every $n\geq1$, let $P_n$ be an
irreducible Markov transition operator on $V_n$, with invariant probability
measure $\widehat\pi_n$. Let $(T_n)_{n\geq1}$ be a strictly increasing
sequence of integer times with $T_1=0$, and define
\[
 L_n=T_{n+1}-T_n.
\]
All these sets, operators and times are deterministic.
We consider the Markov chain started in $V_1$ which uses $P_n$ for the
transitions from $t$ to $t+1$ when $T_n\leq t<T_{n+1}$. Thus the chain
remains in $V_n$ throughout this interval, including at time $T_{n+1}$.
The next transition uses $P_{n+1}$, which is defined on the larger set
$V_{n+1}$. We call these time intervals plateaux.

For every $n\geq1$, let $b_n\geq0$ be an integer, and put
$u_n=(L_n-b_n)_+$. We shall use two estimates on the walk during
each plateau. Assume that there are constants $\alpha>0$ and $B<\infty$,
independent of $n$, such that
\begin{align}
 P_n^m(x,z)&\geq\alpha\widehat\pi_n(z),
 &&x,z\in V_n,\ m\in\N,\ m\geq b_n,\label{eq:plateau-minorization}\\
 \sum_{k=0}^{u}P_n^k(z,z)&\leq
 B(1+u\widehat\pi_n(z)),
 &&z\in V_n,\ u\in\N,\ u\geq1.\label{eq:plateau-green}
\end{align}
The first inequality bounds the transition probabilities from below
after time $b_n$. The second bounds the expected number of visits to
$z$ up to time $u$ when the walk with transition operator $P_n$ starts
from $z$. We use it to control the second moment of the number of visits
within a plateau.

\begin{theorem}[Plateau recurrence]\label{thm:plateau}
Let $(V_n)_{n\geq1}$, $(P_n)_{n\geq1}$ and $(T_n)_{n\geq1}$ be as above.
Assume that~\eqref{eq:plateau-minorization} and~\eqref{eq:plateau-green}
hold with the same constants $\alpha>0$ and $B<\infty$ for every $n\geq1$.
Suppose that
\begin{equation}\label{eq:plateau-divergence}
 \forall z\in V_\infty,\qquad
 \sum_{n:\,z\in V_n}u_n\widehat\pi_n(z)=\infty.
\end{equation}
Then, for every initial distribution on $V_1$, every vertex of
$V_\infty$ is visited infinitely often almost surely.
\end{theorem}

To prove the theorem, we estimate the first two moments of the number of
visits during a plateau and apply the Markov property at the times $T_n$.
This does not require independence between visits in different plateaux.
When the series in~\eqref{eq:plateau-divergence} converges, the theorem
does not say whether the walk is recurrent or transient.

We apply this result to growing boxes. For $R\geq1$, let
\[
 B_R=\{x\in\Z^d:\|x\|_\infty\leq R\}.
\]
Let $(R_n)_{n\geq1}$ be an increasing sequence of integer radii tending
to infinity and take $V_n=B_{R_n}$. For each $n\geq1$, let $c_n$ be
symmetric conductances on the nearest-neighbour edges inside $B_{R_n}$.
Assume that there is a constant $\kappa\in(0,1]$ such that
\begin{equation}\label{eq:box-ellipticity}
 \forall n\geq1,\ \forall e\text{ in }B_{R_n},\qquad
 \kappa\leq c_n(e)\leq\kappa^{-1}.
\end{equation}
Write $\widetilde P_n$ for the conductance walk on the induced finite
graph, and let $P_n=(I+\widetilde P_n)/2$. The sequence $(c_n)$ need
not be monotone.

\begin{corollary}[Growing boxes]\label{cor:boxes}
Let $d\geq3$ and assume~\eqref{eq:box-ellipticity}. There exists a
constant $a=a(d,\kappa)>0$ with the following property. Set
$b_n=\lceil aR_n^2\rceil$ for each $n\geq1$. If
\begin{equation}\label{eq:box-series}
 \sum_{n\geq1}\frac{(L_n-b_n)_+}{R_n^d}=\infty,
\end{equation}
then every vertex of $\Z^d$ is visited infinitely often almost surely.
In particular, if $L_n\geq2b_n$ for all sufficiently large $n$, it is
enough to assume that $\sum_n L_n/R_n^d=\infty$.
\end{corollary}

For example, if $R_n=2^n$ and $L_n=\lceil 2^{dn}/n\rceil$, then
$L_n/b_n\to\infty$ and $L_n/R_n^d\sim1/n$. Hence the series diverges
and Corollary~\ref{cor:boxes} applies. This gives a recurrent walk on
growing boxes whose union is $\Z^d$, although the homogeneous simple
random walk on $\Z^d$ is transient for $d\geq3$.

The proofs of Theorems~\ref{thm:finite-stability} and~\ref{thm:transience}
use an inequality comparing two Dirichlet forms. We first prove it in
finite volume and then pass to infinite volume by exhaustion. The plateau
result uses the two estimates stated above and is proved separately.

\subsection{Historical comments}\label{subsec:history}

Amir, Benjamini, Gurel-Gurevich and Kozma~\cite{ABGGK20} studied random
walks in a changing environment. They stated four conjectures
on recurrence and transience for monotone nonadaptive conductances
\cite[Conjecture~7.1]{ABGGK20}. They proved the corresponding
statements on the half-line and some of them on trees. They also gave
an adaptive example on $\Z^2$ in which the walk is transient although
the conductances are non-decreasing and take values in $\{1,2\}$. This example shows why
one must distinguish an environment fixed in advance from one which
depends on the trajectory of the walk.

Under isoperimetric assumptions, one can also approach transience through
heat-kernel bounds. Dembo, Huang, Morris and Peres used evolving sets to obtain
such bounds and deduce transience~\cite{DHMP17}. Dembo, Huang and Zheng proved
Gaussian estimates for constant-speed continuous-time walks with
non-decreasing reversing measures, assuming uniform ellipticity, volume
regularity and a uniform Poincar\'e inequality~\cite{DHZ19}. Park and
Ray proved recurrence and transience results under boundedness and
summability assumptions on the changes in edge resistance. Their
results also allow adaptive environments~\cite{ParkRay26}.

These questions on recurrence and transience in monotone environments
were at the origin of my PhD thesis. I am happy to see them resolved.
The work leading to this paper began before September 2026.
On 23 September 2026, Li and Park posted the preprint~\cite{LP26}, before
the circulation of the present manuscript. Their Theorems~1.6 and~1.8
prove the two qualitative conclusions stated in this paper.
I acknowledge their priority for these recurrence and transience results.
Their results also cover zero conductances and vertices of infinite
degree, provided the relevant total conductances are finite.
They also prove continuous-time versions and give counterexamples to
the two conjectures for decreasing conductances. Their proofs for
increasing conductances use accumulated electrical currents. In
addition to recurrence and transience, they give quantitative bounds
on escape probabilities and expected occupation times.

Our proofs compare the hitting probabilities of the time-inhomogeneous
walk with those of the walk whose conductances are fixed at the starting
time. In Theorem~\ref{thm:finite-stability} we bound the absolute
difference of the two exit probabilities in a finite set. In infinite
volume, Inequality~\eqref{eq:transience-stability} gives a lower bound
on the probability of avoiding $o$ forever. To obtain these estimates,
we compare two Dirichlet forms and sum the changes of the equilibrium
potentials. In both bounds the capacity enters through its remaining
increase after the starting time. I thank Pierre Mathieu for pointing out that the argument
extends to a more general setting.

For growing Euclidean domains, Dembo, Huang and Sidoravicius obtained
criteria involving the integral of the inverse volume~\cite{DHS14}.
Their discrete-time recurrence theorem requires additional assumptions
on the shape of the domains and separation-of-scales assumptions on
the growth. Huang proved transition-probability estimates and
sufficient criteria for transience and weak recurrence on growing
graphs~\cite{Huang19}.

Kumamoto, Kijima and Shirai~\cite{KKS24} also consider walks that spend
a prescribed number $L_n$ of steps on each graph. For their reflected
walk of period two on the boxes $[-n,n]^d\cap\Z^d$, $d\geq4$, the expected
number of returns to the origin is infinite if and only if
$\sum_n L_n/n^d$ diverges (Theorem~5.2 of the full version). On growing
hypercubes, they also prove almost-sure repeated visits to each fixed
vertex under the assumption $L_n\geq c n2^n$ for all large $n$, with
$c>0$ \cite[Proposition~22]{KKS24}.

In Theorem~\ref{thm:plateau}, we assume the
two finite-volume estimates~\eqref{eq:plateau-minorization}
and~\eqref{eq:plateau-green}. We deduce that every vertex is visited
infinitely often almost surely when the series
in~\eqref{eq:plateau-divergence} diverges.

In Corollary~\ref{cor:boxes}, the walk on each box is lazy, the radii
may increase at arbitrary rates, and the uniformly elliptic conductances
may vary from one plateau to the next. We verify the two estimates on
each box using finite-volume Nash and Poincar\'e inequalities. This gives
a sufficient condition for almost-sure recurrence in dimensions $d\geq3$.
The relation between Gaussian estimates, volume regularity and
Poincar\'e inequalities was studied by Delmotte~\cite{Delmotte99}.
Quantitative merging estimates obtained from functional inequalities
for non-decreasing invariant measures are given in~\cite{Moumeni26}.
We include the estimates needed for the plateau theorem in
Section~\ref{sec:plateaux}.

In Section~\ref{sec:examples} we apply the estimates to walks on fixed
graphs, including the square lattice. In Section~\ref{sec:twoforms},
we prove the comparison estimate and deduce recurrence. Section~\ref{sec:transience}
is devoted to transience. In Section~\ref{sec:plateaux}, we prove the
plateau theorem and the estimates for growing boxes. In
Appendix~\ref{sec:electrical}, I give my original proof of the two-form
inequality, based on an interpolation between the two electrical networks.

\section{Examples and consequences}\label{sec:examples}

We first consider examples in which the graph is fixed and its
conductances increase.
We will consider growing boxes in Section~\ref{sec:plateaux}. There,
the domains increase, while the transition operator remains fixed
during each plateau.

\subsection{Uniformly elliptic recurrent graphs}

Let $G$ be a graph which is recurrent for unit conductances. Assume
that there are constants $c_->0$ and $c_+<\infty$ such that, for
every edge $e$ and every $t\geq0$,
\[
 0<c_-\leq c_t(e)\leq c_{t+1}(e)\leq c_+<\infty.
\]
Set $\bar c(e)=c_+$ for every edge. Multiplication of all
conductances by the same constant preserves recurrence, so
$(G,\bar c)$ is recurrent. Theorem~\ref{thm:recurrence} therefore
applies. In particular, the walk visits every vertex infinitely often
almost surely when $G=\Z$ or $G=\Z^2$.

Only the upper bound $c_+$ was needed to obtain a recurrent upper
network. Theorem~\ref{thm:recurrence} also allows upper conductances
which depend on the edge. We will use the lower bound $c_-$ in the
following estimate.

Let us give the corresponding escape estimate on $\Z^2$. For
$R\geq4$, set $B_R=\{z\in\Z^2:\|z\|_\infty\leq R\}$.
Comparing the energies with those of the unit lattice gives, for
every $s\geq0$ and every $x\in B_{R/2}\setminus\{0\}$,
\[
 R_s^{\mathrm{free}}(x,0)\leq C\log(1+\|x\|_\infty),\qquad
 \Reff^{\bar c}(0,B_R^c)\geq c\log R.
\]
In the first bound, $c_s(e)\geq c_-$ controls the resistance between
$x$ and $0$. In the second, $\bar c\equiv c_+$ allows us to use the
resistance from $0$ to the exterior of the box in the unit network.
It follows from Corollary~\ref{cor:escape} that
\begin{equation}\label{eq:planar-example}
 \PP_{s,x}(\tau_{B_R^c}^s<\tau_0^s)
 \leq C\sqrt{\frac{\log(1+\|x\|_\infty)}{\log R}}.
\end{equation}
The constants depend only on $c_-$ and $c_+$ and are therefore
uniform in $s$ and in the sequence of conductances. The
resistance estimates for the unit lattice can be found in
\cite[Chapter~2]{LP16book}. Inequality~\eqref{eq:planar-example}
contains a square root, which comes from the energy estimate.
For simple random walk on $\Z^2$, the escape probability itself is
of order $\log(1+\|x\|_\infty)/\log R$. This follows by stopping the
planar potential kernel when the walk first hits $0$ or exits $B_R$
\cite[Chapter~4]{LP16book}.

\subsection{A transient initial network}

Let $G=\Z^d$, where $d\geq3$. Assume that, for every edge $e$ and
every $t\geq0$,
\[
 1\leq c_t(e)\leq c_{t+1}(e)\leq\bar c(e)<\infty.
\]
Since $c_0$ dominates the unit conductances on $\Z^d$, Rayleigh
monotonicity gives transience of $(G,c_0)$. Thus
Theorem~\ref{thm:transience} applies, and the expected number of
visits to every fixed vertex is finite. Moreover,
Inequality~\eqref{eq:last-exit} gives a strictly positive lower bound
on $\liminf_{t\to\infty}q_t(o)$, where $q_t(o)$ is the probability
that a walk started at $o$ at time $t$ never returns. The upper
conductance bound may depend on the edge. There is no further restriction on the times of the
increases or on their sizes.

The same argument applies to any locally finite graph for which
$(G,c_0)$ is transient. Indeed, since each vertex has finitely many
neighbours, the condition $c_\infty(e)<\infty$ on every edge gives
a finite upper bound on the total conductance at that vertex.
We use this local bound in the proof when summing escape
probabilities over the neighbours of the target vertex.

\subsection{The role of the capacity increment}

Assume that the conductances are constant after a deterministic
time $T$. Fix a finite domain $D$. For every $s\geq T$, we have
$C_*^D-C_s^D=0$, and Theorem~\ref{thm:finite-stability} gives
$p_s^D=u_s^D$. For $s<T$, the same theorem bounds the difference
in terms of the remaining increase of capacity and the dual energy
at time $s$. In particular, we do not need to estimate the
individual increments.

Let us see why the bound can be strict. Take
$c_t=a_t c$, where $a_t$ is a positive non-decreasing scalar sequence
with finite limit $a_\infty$. Then
\[
 K_t(x,y)=\frac{a_t c(x,y)}{\sum_z a_t c(x,z)}=K_c(x,y).
\]
Thus the transition operator does not depend on time, and
$p_s^D=u_s^D$ for every $s$. On the other hand, $C_s^D=a_s C_c^D$, where $C_c^D$ is the capacity
for the conductances $c$. If $a_s<a_\infty$ and $x\in F=D\setminus\{o\}$ has a neighbour
in $F$, then $\delta_x K_s^F\ne0$,
so its dual energy is positive. The first bound in
\eqref{eq:finite-stability} is then strictly positive even though
the two exit probabilities agree.

\subsection{Why the environment is prescribed in advance}

Suppose now that the conductances are random and satisfy
Equality~\eqref{eq:quenched}. After conditioning on the full
environment, the transition from time $t$ to time $t+1$ is still
given by $K_t$. The deterministic results therefore apply to the
conditional law of the walk.

For an adaptive rule, the future conductances may reveal which
edge the walk has traversed, and Equality~\eqref{eq:quenched} may
fail. In particular, the conditional law of the killed walk need
not be given by the product of the restricted kernels $K_t^F$.
This product is what allows us in Section~\ref{sec:twoforms} to
sum the changes of the frozen potentials against the law of the
killed walk.

In the adaptive example of Amir, Benjamini, Gurel-Gurevich and
Kozma~\cite{ABGGK20}, the conductances on $\Z^2$ are non-decreasing
and take values in $\{1,2\}$, yet the walk is transient.
We thus keep the deterministic assumption, or the conditional
transition identity~\eqref{eq:quenched} when the environment is
random.

\section{Two-form comparison and finite-volume estimates}\label{sec:twoforms}

We first compare the two exit probabilities in Theorem~\ref{thm:finite-stability}. We then use this comparison to deduce the escape bound of Corollary~\ref{cor:escape} and the recurrence statement of Theorem~\ref{thm:recurrence}. We consider deterministic conductances satisfying
\[
0<c_t(e)\leq c_{t+1}(e)\leq\bar c(e)<\infty,
\qquad e\in E,\quad t\geq0.
\]
Let $D$ be a finite connected set containing $o$. Set $F=D\setminus\{o\}$ and $A=\{o\}\cup D^c$. We take $D$ to contain at least one vertex other than $o$, since the starting point for the exit estimate lies in $F$. We extend test functions on $F$ by zero on $A$. For a real signed measure $\nu$ on $F$ and a function $f$, we write $\nu(f)=\sum_{y\in F}\nu(y)f(y)$.

\subsection{Potentials and dual energies}

For every $t\geq0$, recall that $u_t^D$ is the unique function equal to $0$ at $o$, equal to $1$ on $D^c$, and $c_t$-harmonic on $F$. This function minimizes $\cE_t(f,f)$ among functions with the same boundary values. Thus,
\begin{equation}\label{eq:twoforms-orthogonality}
\cE_t(u_t^D,f)=0
\qquad\text{whenever }f|_A=0,
\end{equation}
and
\[
C_t^D=\cE_t(u_t^D,u_t^D)
=\frac{1}{\Reff^{c_t}(o,D^c)}.
\]
The conductances are non-decreasing, so the Dirichlet forms are ordered. The same order holds for their infima over functions with the prescribed boundary values. Hence $C_t^D\leq C_{t+1}^D\leq C_{\bar c}^D$. Consequently, the limit $C_*^D=\lim_{t\to\infty}C_t^D$ exists and is finite.

For every signed measure $\nu$ on $F$, let
\begin{equation}\label{eq:twoforms-dual-definition}
H_t^D(\nu)
=\sup_{f|_A=0}\bigl\{2\nu(f)-\cE_t(f,f)\bigr\}.
\end{equation}
The form $\cE_t$ is positive definite on the functions vanishing on $A$. Indeed, a function of zero energy is constant along every path in $D$, and it vanishes at $o$. Hence there exists a unique function $\phi_t^\nu$, vanishing on $A$, such that
\begin{equation}\label{eq:twoforms-dual-equation}
\cE_t(\phi_t^\nu,f)=\nu(f)
\qquad\text{for all }f|_A=0.
\end{equation}
Completing the square in~\eqref{eq:twoforms-dual-definition}, we find
\begin{equation}\label{eq:twoforms-dual-value}
H_t^D(\nu)=\nu(\phi_t^\nu)
=\cE_t(\phi_t^\nu,\phi_t^\nu).
\end{equation}
Moreover, the definition gives $H_t^D(\nu)\geq H_{t+1}^D(\nu)$, because $\cE_t\leq\cE_{t+1}$ and the supremum is taken over the same set of functions.

\begin{lemma}[Two-form comparison]\label{lem:twoform-comparison}
For every $t\geq0$ and every signed measure $\nu$ on $F$,
\begin{equation}\label{eq:twoforms-comparison}
\bigl|\nu(u_{t+1}^D-u_t^D)\bigr|^2
\leq
(C_{t+1}^D-C_t^D)
\bigl(H_t^D(\nu)-H_{t+1}^D(\nu)\bigr).
\end{equation}
\end{lemma}

\begin{proof}
For $\lambda\in\R$, let
\[
J_t(\lambda)
=\inf_{\substack{f(o)=0\\f|_{D^c}=1}}
\bigl\{\cE_t(f,f)-2\lambda\nu(f)\bigr\}.
\]
Write an admissible function as $u_t^D+g$, where $g|_A=0$. By Equalities~\eqref{eq:twoforms-orthogonality} and~\eqref{eq:twoforms-dual-equation}, we have
\begin{align*}
&\cE_t(u_t^D+g,u_t^D+g)-2\lambda\nu(u_t^D+g)\\
&\qquad=C_t^D+\cE_t(g,g)-2\lambda\nu(u_t^D)-2\lambda\nu(g)\\
&\qquad=C_t^D-2\lambda\nu(u_t^D)-\lambda^2H_t^D(\nu)
  +\cE_t(g-\lambda\phi_t^\nu,g-\lambda\phi_t^\nu).
\end{align*}
The last term is nonnegative and is zero when $g=\lambda\phi_t^\nu$. Thus the minimizing function is $u_t^D+\lambda\phi_t^\nu$, and we find
\[
J_t(\lambda)
=C_t^D-2\lambda\nu(u_t^D)-\lambda^2H_t^D(\nu).
\]
For every real $\lambda$, we also have $J_t(\lambda)\leq J_{t+1}(\lambda)$, by $\cE_t\leq\cE_{t+1}$. Subtracting the two expressions gives
\[
(C_{t+1}^D-C_t^D)
-2\lambda\nu(u_{t+1}^D-u_t^D)
+\lambda^2\bigl(H_t^D(\nu)-H_{t+1}^D(\nu)\bigr)
\geq0.
\]
This polynomial is nonnegative for every real $\lambda$. If its quadratic coefficient is positive, Inequality~\eqref{eq:twoforms-comparison} follows from the sign of the discriminant. If that coefficient is zero, the linear coefficient must vanish as well, since $\lambda$ may have either sign.
\end{proof}

For every $t\geq0$, let $K_t^F$ be the restriction of $K_t$ to $F\times F$. It is the transition matrix of the walk killed upon entering $A$. The next lemma shows that the dual energy decreases under this transition. It holds without a laziness assumption.

\begin{lemma}[Contraction of the dual energy]\label{lem:killed-dual-contraction}
For every $t\geq0$ and every signed measure $\nu$ on $F$,
\begin{equation}\label{eq:twoforms-contraction}
H_t^D(\nu K_t^F)\leq H_t^D(\nu).
\end{equation}
Moreover, for $x\in F$,
\begin{equation}\label{eq:twoforms-resistance-comparison}
H_t^D(\delta_x)=R_t^D(x,A)
\leq R_t^{\mathrm{free}}(x,o).
\end{equation}
\end{lemma}

\begin{proof}
Fix $t\geq0$ and set $P=K_t^F$. We work in $\ell^2(F,\pi_t)$. By reversibility, $P$ is self-adjoint, and for functions extended by zero on $A$ we have
\[
\cE_t(f,f)=\langle f,(I-P)f\rangle_{\pi_t}
\qquad\text{for }f|_A=0.
\]
Every vertex of $F$ can be joined to $o$ by a path in $D$ with positive transition probabilities. Since $F$ is finite, the lengths of these paths have a finite maximum and their probabilities have a positive minimum. Let $m$ be the maximum of these lengths and $a$ the minimum of these probabilities. Whatever the starting vertex, the frozen walk reaches $A$ within $m$ steps with probability at least $a$. Thus every row sum of $P^m$ is at most $1-a$, and the spectral radius of $P$ is strictly less than one.

For every $y\in F$, let $g(y)=\nu(y)/\pi_t(y)$. By reversibility, $Pg$ is the density of $\nu P$ with respect to $\pi_t$. Using Equalities~\eqref{eq:twoforms-dual-equation} and~\eqref{eq:twoforms-dual-value}, we obtain
\[
H_t^D(\nu)=\langle g,(I-P)^{-1}g\rangle_{\pi_t},
\qquad
H_t^D(\nu P)=\langle g,P^2(I-P)^{-1}g\rangle_{\pi_t}.
\]
For every $\theta$ in the spectrum of $P$, we have $\theta\in(-1,1)$ and $\theta^2/(1-\theta)\leq1/(1-\theta)$. Applying the spectral theorem to the two expressions above gives Inequality~\eqref{eq:twoforms-contraction}, including when $P$ has negative eigenvalues.

To prove the resistance comparison, we return to Equality~\eqref{eq:twoforms-dual-definition}. We first optimize over scalar multiples of a fixed nonzero function $f$, then take the supremum over $f$. We obtain
\[
H_t^D(\delta_x)
=\sup_{\substack{f|_A=0\\f\ne0}}
\frac{f(x)^2}{\cE_t(f,f)}.
\]
The last supremum is the resistance between $x$ and the grounded set $A$. Moreover, every function in this supremum has finite energy and satisfies $f(o)=0$, so it is also admissible in the definition of $R_t^{\mathrm{free}}(x,o)$. This proves the second assertion.
\end{proof}

\subsection{The killed evolution}

For every $s\geq0$ and every $t\geq s$, set
\[
K_{s,t}^F=K_s^F K_{s+1}^F\cdots K_{t-1}^F,
\qquad K_{s,s}^F=I.
\]
We first prove that the evolving walk reaches $A$ almost surely. The probability that it survives until time $T$ will bound the remainder in the comparison formula.

\begin{lemma}[Finite-volume absorption]\label{lem:finite-volume-absorption}
For each $s\geq0$, there are an integer $m\geq1$ and a number $a>0$, depending on $s$ and $D$, such that
\begin{equation}\label{eq:twoforms-absorption}
\sup_{y\in F}\delta_y K_{s,s+km}^F(F)
\leq(1-a)^k,
\qquad k\geq0.
\end{equation}
In particular, $\tau_A^s<\infty$ almost surely under $\PP_{s,y}$ for every $y\in F$.
\end{lemma}

\begin{proof}
Fix $s\geq0$. For each $y\in F$, choose a path from $y$ to $o$ contained in $D$, and let $m$ be the largest length of these paths. By monotonicity and the upper bound on the conductances, for every edge $(z,w)$ on one of the chosen paths and every $t\geq s$,
\[
K_t(z,w)=\frac{c_t(z,w)}{\pi_t(z)}
\geq\frac{c_s(z,w)}{\pi_{\bar c}(z)}>0.
\]
Here $\pi_{\bar c}(z)=\sum_{w\sim z}\bar c(z,w)$ is finite by local finiteness. For each chosen path, multiply the lower bounds along its edges and choose $0<a<1$ no larger than the smallest of these products. Thus, for every starting time $t\geq s$ and every starting point in $F$, the probability of hitting $o$ within $m$ steps is at least $a$. Applying the Markov property at successive intervals of length $m$ gives Inequality~\eqref{eq:twoforms-absorption}.
\end{proof}

Now fix $s\geq0$ and $x\in F$. For every $t\geq s$, let
\[
\mu_t=\delta_xK_{s,t}^F.
\]
Thus $\mu_t$ is the law of $X_t$ restricted to the event $\{\tau_A^s>t\}$, and $\mu_{t+1}=\mu_tK_t^F$. We recall that
\[
p_t^D(y)=\PP_{t,y}(\tau_{D^c}^t<\tau_{\{o\}}^t).
\]
For every $y\in F$, write $b_t(y)=K_t(y,D^c)$. Conditioning on the first step of the evolving walk and of the walk frozen at time $t$, respectively, gives
\[
p_t^D=b_t+K_t^Fp_{t+1}^D,
\qquad
u_t^D=b_t+K_t^Fu_t^D.
\]
Subtracting the second equality from the first gives
\[
p_t^D-u_t^D
=K_t^F(p_{t+1}^D-u_{t+1}^D)
 +K_t^F(u_{t+1}^D-u_t^D).
\]
Apply $\mu_t$ to both sides. Since $\mu_tK_t^F=\mu_{t+1}$, the first term on the right is $\mu_{t+1}(p_{t+1}^D-u_{t+1}^D)$. The sum from $t=s$ to $T-1$ telescopes. Thus, for every $T>s$,
\begin{equation}\label{eq:twoforms-duhamel}
p_s^D(x)-u_s^D(x)
=\sum_{t=s}^{T-1}\mu_{t+1}(u_{t+1}^D-u_t^D)
+\mu_T(p_T^D-u_T^D).
\end{equation}
Both $p_T^D$ and $u_T^D$ take values in $[0,1]$. Hence the absolute value of the last term is at most $\mu_T(F)$, which tends to zero by Lemma~\ref{lem:finite-volume-absorption}.

\begin{proof}[Proof of Theorem~\ref{thm:finite-stability}]
For every $t\geq s$, apply Lemma~\ref{lem:twoform-comparison} with $\nu=\mu_{t+1}$, and set
\[
\Delta_t
=H_t^D(\mu_{t+1})-H_{t+1}^D(\mu_{t+1})\geq0.
\]
Summing the resulting bounds and applying the Cauchy--Schwarz inequality, we find
\begin{equation}\label{eq:twoforms-sum-bound}
\sum_{t=s}^{T-1}\bigl|\mu_{t+1}(u_{t+1}^D-u_t^D)\bigr|
\leq
\left((C_T^D-C_s^D)\sum_{t=s}^{T-1}\Delta_t\right)^{1/2}.
\end{equation}
Let us apply Lemma~\ref{lem:killed-dual-contraction} at time $t+1$. We have
\[
H_{t+1}^D(\mu_{t+2})
=H_{t+1}^D(\mu_{t+1}K_{t+1}^F)
\leq H_{t+1}^D(\mu_{t+1}).
\]
We can therefore replace $H_{t+1}^D(\mu_{t+1})$ in the subtracted term of $\Delta_t$ by $H_{t+1}^D(\mu_{t+2})$ to obtain an upper bound. The sum now telescopes.
\begin{align*}
\sum_{t=s}^{T-1}\Delta_t
&\leq\sum_{t=s}^{T-1}
\bigl(H_t^D(\mu_{t+1})-H_{t+1}^D(\mu_{t+2})\bigr)\\
&=H_s^D(\mu_{s+1})-H_T^D(\mu_{T+1})\\
&\leq H_s^D(\delta_xK_s^F).
\end{align*}
Together with Inequality~\eqref{eq:twoforms-sum-bound}, this proves that the series in Equality~\eqref{eq:twoforms-duhamel} converges absolutely. We let $T$ tend to infinity in~\eqref{eq:twoforms-sum-bound} and use the convergence to zero of the remainder in~\eqref{eq:twoforms-duhamel}. This gives
\begin{equation}\label{eq:twoforms-fine-stability}
\bigl|p_s^D(x)-u_s^D(x)\bigr|
\leq
\sqrt{(C_*^D-C_s^D)H_s^D(\delta_xK_s^F)}.
\end{equation}
Finally, recall that $C_*^D\leq C_{\bar c}^D$. Moreover,
\[
H_s^D(\delta_xK_s^F)
\leq H_s^D(\delta_x)=R_s^D(x,A),
\]
by Lemma~\ref{lem:killed-dual-contraction}. Substituting these two bounds in the preceding inequality gives
\begin{equation}\label{eq:twoforms-grounded-stability}
\bigl|p_s^D(x)-u_s^D(x)\bigr|
\leq
\sqrt{(C_{\bar c}^D-C_s^D)R_s^D(x,A)}.
\end{equation}
This proves the theorem.
\end{proof}

\subsection{Escape probabilities and recurrence}

\begin{proof}[Proof of Corollary~\ref{cor:escape}]
The potential $u_s^D$ has finite energy and satisfies $u_s^D(o)=0$. We may therefore use it in the definition of the free resistance, which gives
\[
u_s^D(x)^2
\leq R_s^{\mathrm{free}}(x,o)\cE_s(u_s^D,u_s^D)
=R_s^{\mathrm{free}}(x,o)C_s^D
\leq R_s^{\mathrm{free}}(x,o)C_{\bar c}^D.
\]
Combining this bound with Inequalities~\eqref{eq:twoforms-grounded-stability} and~\eqref{eq:twoforms-resistance-comparison}, we obtain
\begin{align*}
\PP_{s,x}(\tau_{D^c}^s<\tau_{\{o\}}^s)
&=p_s^D(x)\\
&\leq 2\sqrt{R_s^{\mathrm{free}}(x,o)C_{\bar c}^D}
=2\sqrt{\frac{R_s^{\mathrm{free}}(x,o)}{\Reff^{\bar c}(o,D^c)}}.
\end{align*}
The numerator is bounded by the sum of the edge resistances along any fixed path from $x$ to $o$ and hence is finite, independently of $D$.
\end{proof}

\begin{proof}[Proof of Theorem~\ref{thm:recurrence}]
Fix a deterministic time $s\geq0$ and vertices $x,o\in V$ with $x\ne o$. Let $(D_n)$ be an increasing exhaustion of $V$ by finite connected sets containing $o$ and $x$. Since $(G,\bar c)$ is recurrent, we have
\[
\Reff^{\bar c}(o,D_n^c)\longrightarrow\infty.
\]
For every $n$, Lemma~\ref{lem:finite-volume-absorption} shows that the chain almost surely hits $o$ or exits $D_n$ in finite time. In particular, on the event that it never hits $o$, it must exit $D_n$ before hitting $o$. Corollary~\ref{cor:escape} then gives
\begin{align*}
\PP_{s,x}(\tau_{\{o\}}^s=\infty)
&\leq\PP_{s,x}(\tau_{D_n^c}^s<\tau_{\{o\}}^s)\\
&\leq2\sqrt{\frac{R_s^{\mathrm{free}}(x,o)}{\Reff^{\bar c}(o,D_n^c)}}
\longrightarrow0.
\end{align*}
Thus the walk hits $o$ almost surely from every $x\ne o$ and every deterministic starting time $s$. For $x=o$, the hitting time is $s$ by definition.

We now deduce that the walk visits the target infinitely often. Let $\sigma_0$ be the first visit to $o$ at or after the prescribed starting time, and define recursively
\[
\sigma_{k+1}=\inf\{t\geq\sigma_k+1:X_t=o\}.
\]
We have just proved that $\sigma_0<\infty$ almost surely. On $\{\sigma_k<\infty\}$, condition on the time and position at $\sigma_k+1$. By the strong Markov property of the space-time chain, the hitting result applies from that point and time. Hence $\sigma_{k+1}<\infty$ almost surely on $\{\sigma_k<\infty\}$. By induction, the walk visits $o$ infinitely often almost surely. Since $V$ is countable, this holds simultaneously for every vertex.
\end{proof}

\section{Transience and escape from a vertex}\label{sec:transience}

In this section, we prove Theorem~\ref{thm:transience}. We assume
that the initial network $(G,c_0)$ is transient and that the
conductances satisfy~\eqref{eq:monotone}. By Rayleigh monotonicity,
$(G,c_t)$ is then transient for every $t\geq0$.

Let $o\in V$ be fixed. We use the functions $u_t$, the capacities
$\kappa_t(o)$ and the dual energies $H_t^o$ defined in the
introduction. Since $1-u_t$ is the equilibrium potential of
$\{o\}$ in the network $(G,c_t)$, we have
\begin{equation}\label{eq:infinite-capacity}
 \cE_t(u_t,u_t)=\kappa_t(o),\qquad
 0<\kappa_0(o)\leq\kappa_t(o)\leq\pi_\infty(o).
\end{equation}
By monotonicity, $\kappa_t(o)$ increases to its finite positive
limit $\kappa_*(o)$.

\subsection{The two-form inequality in infinite volume}

We first check that $H_t^o(\nu)$ is finite for every finitely
supported signed measure $\nu$. For each $x$ in the support of
$\nu$, choose a simple path $\gamma_x$ from $x$ to $o$. If
$f(o)=0$, write $f(x)$ as the sum of its differences along this
path. The Cauchy--Schwarz inequality gives
\[
 |f(x)|\leq
 \left(\sum_{e\in\gamma_x}\frac1{c_t(e)}\right)^{1/2}
 \cE_t(f,f)^{1/2}.
\]
Multiplying by $|\nu(x)|$ and summing over the support of
$\nu$, we bound the linear term in the definition of $H_t^o(\nu)$
by a constant times the square root of the energy. Optimizing
over the energy gives
\begin{equation}\label{eq:infinite-dual-finite}
 H_t^o(\nu)\leq
 \left[\sum_x|\nu(x)|
 \left(\sum_{e\in\gamma_x}\frac1{c_t(e)}\right)^{1/2}\right]^2.
\end{equation}
In particular, if $y$ is a neighbour of $o$, we may take the path
consisting of the edge $(o,y)$. We then obtain
\begin{equation}\label{eq:neighbour-dual}
 H_t^o(\delta_y)\leq\frac1{c_t(o,y)}\leq\frac1{c_0(o,y)}.
\end{equation}

\begin{lemma}[Infinite-volume comparison]\label{lem:infinite-twoforms}
For every $t\geq0$ and every finitely supported signed measure $\nu$,
\begin{equation}\label{eq:infinite-twoforms}
 |\nu(u_{t+1}-u_t)|^2
 \leq(\kappa_{t+1}(o)-\kappa_t(o))
       (H_t^o(\nu)-H_{t+1}^o(\nu)).
\end{equation}
\end{lemma}

\begin{proof}
The functions $u_t$, $u_{t+1}$ and the test functions defining the
dual energies vanish at $o$. Thus the mass of $\nu$ at $o$ does
not contribute to either side, and we may assume $\nu(o)=0$.
Let $D_n\uparrow V$ be an exhaustion by finite connected sets
containing $o$ and $\supp(\nu)$. For each $n$, apply
Lemma~\ref{lem:twoform-comparison} in $D_n$.

Let us pass to the limit in $n$, keeping $t$ fixed. For the
network frozen at time $t$, we have, for each $x\in V$,
\[
 u_t^{D_n}(x)\downarrow u_t(x),\qquad
 C_t^{D_n}\downarrow\kappa_t(o).
\]
Indeed, the frozen walk leaves every finite set almost surely.
Hence the events of exiting $D_n$ before hitting $o$ decrease, up to
a null set, to the event of never hitting $o$. The second limit is
the exhaustion formula for capacity \cite[Chapter~2]{LP16book}.
For the dual energies, we have
\[
 H_t^{D_n}(\nu)\uparrow H_t^o(\nu).
\]
Indeed, the test spaces in the finite-volume supremum increase,
and every compactly supported function vanishing at $o$ belongs
to one of them. Inequality~\eqref{eq:infinite-dual-finite} bounds this supremum.
Taking these limits at both $t$ and $t+1$ in the finite-volume
inequality proves Inequality~\eqref{eq:infinite-twoforms}.
\end{proof}

\begin{remark}[The boundary in the dual energy]
For every $x\ne o$, the same exhaustion gives
\[
 H_t^o(\delta_x)
 =\lim_{n\to\infty}
 \Reff^{c_t}\bigl(x,\{o\}\cup D_n^c\bigr).
\]
In this expression, the exterior of $D_n$ is grounded together
with $o$. For a transient network, this limit may differ from the
free resistance
$R_t^{\mathrm{free}}(x,o)$ defined in~\eqref{eq:free-resistance}.
For neighbours of $o$, we will use
Inequality~\eqref{eq:neighbour-dual}, obtained from the single edge
joining the two vertices.
\end{remark}

For every $t\geq0$, let $K_t^{(o)}$ be the kernel $K_t$ killed
upon hitting $o$. As in finite volume, we extend functions on
$V\setminus\{o\}$ by zero at $o$.

\begin{lemma}[Killed contraction]\label{lem:infinite-contraction}
For every finitely supported signed measure $\nu$ on $V\setminus\{o\}$,
\begin{equation}\label{eq:infinite-contraction}
 H_t^o(\nu K_t^{(o)})\leq H_t^o(\nu).
\end{equation}
\end{lemma}

\begin{proof}
By reversibility, the killed kernel is a self-adjoint contraction on
$\ell^2(V\setminus\{o\},\pi_t)$, with spectrum contained in
$[-1,1]$. For every compactly supported $f$ vanishing at $o$,
\[
 \cE_t(f,f)=\langle f,(I-K_t^{(o)})f\rangle_{\pi_t}.
\]
For every $\lambda\in[-1,1]$, we have
$\lambda^2(1-\lambda)\leq1-\lambda$. Applying the spectral theorem
to the preceding equality gives
\[
 \cE_t(K_t^{(o)}f,K_t^{(o)}f)\leq\cE_t(f,f).
\]
Moreover, $K_t^{(o)}f$ has finite support, since the graph is
locally finite. It is therefore an admissible test function in
the definition of the dual energy. Using the preceding energy
bound, we find
\begin{align*}
 H_t^o(\nu K_t^{(o)})
 &=\sup_f\{2\nu(K_t^{(o)}f)-\cE_t(f,f)\}\\
 &\leq\sup_f\{2\nu(K_t^{(o)}f)
                    -\cE_t(K_t^{(o)}f,K_t^{(o)}f)\}\\
 &\leq H_t^o(\nu).
\end{align*}
The scalar inequality holds on all of $[-1,1]$, so the contraction
also applies to walks that are not lazy.
\end{proof}

\subsection{Comparison with the frozen escape probability}

Let $s\geq0$ and $x\ne o$. Under $\PP_{s,x}$, write
$\tau=\tau_o^s$. For every $t\geq s$ and every vertex $z$, define
\[
 \mu_t(z)=\PP_{s,x}(X_t=z,\tau>t).
\]
The measure $\mu_t$ has finite support. Indeed, local finiteness
implies that the walk can reach only finitely many vertices in a
fixed number of steps. Moreover, $\mu_s=\delta_x$ and
$\mu_{t+1}=\mu_tK_t^{(o)}$.

For every $t\geq s$, the function $u_t$ is $K_t$-harmonic off $o$
and vanishes at $o$. Conditioning on the next step therefore gives
\[
 \EE_{s,x}[\1_{\{\tau>t+1\}}u_t(X_{t+1})\mid\cF_t]
 =\1_{\{\tau>t\}}u_t(X_t).
\]
Taking expectations, we can now replace $u_t$ by $u_{t+1}$ at
time $t+1$ and keep track of the change in potential. We obtain
\begin{align*}
&\EE_{s,x}[\1_{\{\tau>t+1\}}u_{t+1}(X_{t+1})]
 -\EE_{s,x}[\1_{\{\tau>t\}}u_t(X_t)]\\
&\qquad=\EE_{s,x}\bigl[\1_{\{\tau>t+1\}}
          (u_{t+1}-u_t)(X_{t+1})\bigr]
 =\mu_{t+1}(u_{t+1}-u_t).
\end{align*}
When we sum from $t=s$ to $T-1$, the left-hand side telescopes.
Since $\mu_s=\delta_x$, we have, for every $T>s$,
\begin{equation}\label{eq:infinite-martingale}
 \EE_{s,x}[\1_{\{\tau>T\}}u_T(X_T)]
 =u_s(x)+\sum_{t=s}^{T-1}\mu_{t+1}(u_{t+1}-u_t).
\end{equation}
Apply Lemma~\ref{lem:infinite-twoforms} to $\mu_{t+1}$ at each
time $t$. Summing the resulting bounds and applying the
Cauchy--Schwarz inequality gives
\begin{align*}
 \sum_{t=s}^{T-1}|\mu_{t+1}(u_{t+1}-u_t)|
 &\leq\sqrt{\kappa_T(o)-\kappa_s(o)}\\
 &\quad\times\left[\sum_{t=s}^{T-1}
       (H_t^o(\mu_{t+1})-H_{t+1}^o(\mu_{t+1}))\right]^{1/2}.
\end{align*}
Since $\mu_{t+1}=\mu_tK_t^{(o)}$,
Lemma~\ref{lem:infinite-contraction} gives
$H_t^o(\mu_{t+1})\leq H_t^o(\mu_t)$. We use this bound
for the first term of each increment and obtain
\begin{align*}
 \sum_{t=s}^{T-1}
       (H_t^o(\mu_{t+1})-H_{t+1}^o(\mu_{t+1}))
 &\leq\sum_{t=s}^{T-1}
       (H_t^o(\mu_t)-H_{t+1}^o(\mu_{t+1}))\\
 &=H_s^o(\delta_x)-H_T^o(\mu_T)\\
 &\leq H_s^o(\delta_x).
\end{align*}
We use these two estimates in Equality~\eqref{eq:infinite-martingale}.
Since $0\leq u_T\leq1$, the expectation on the left is at most
the probability that the walk avoids $o$ up to time $T$. Thus,
\[
 \PP_{s,x}(\tau>T)
 \geq u_s(x)-\sqrt{(\kappa_T(o)-\kappa_s(o))H_s^o(\delta_x)}.
\]
Letting $T$ tend to infinity gives
Inequality~\eqref{eq:transience-stability}, since the events
$\{\tau>T\}$ decrease to $\{\tau=\infty\}$. Here we use only that
the frozen capacities increase to $\kappa_*(o)$, without identifying
this limit with the capacity of the limiting network.

\subsection{A positive probability of last departure}

We now prove that the expected number of visits to $o$ is finite.
For every $s\geq0$, let
\[
 \eps_s(o)=\max_{y\sim o}
 \sqrt{\frac{\kappa_*(o)-\kappa_s(o)}{c_0(o,y)}}.
\]
The maximum is taken over finitely many neighbours of $o$,
and each denominator is strictly positive. Since $\kappa_s(o)$ converges to
$\kappa_*(o)$, we have $\eps_s(o)\to0$.
Inequalities~\eqref{eq:transience-stability}
and~\eqref{eq:neighbour-dual} therefore give
\[
 \PP_{s,y}(\tau_o^s=\infty)\geq u_s(y)-\eps_s(o),
 \qquad y\sim o.
\]
Start the walk from $o$ at time $t$ and condition on its first
step. By the definition of $q_t(o)$, the preceding bound at
time $t+1$ gives
\begin{align}
 q_t(o)&=\sum_{y\sim o}K_t(o,y)
            \PP_{t+1,y}(\tau_o^{t+1}=\infty)\nonumber\\
 &\geq\sum_{y\sim o}K_t(o,y)u_{t+1}(y)-\eps_{t+1}(o).
 \label{eq:last-departure-average}
\end{align}
For each neighbour $y$ of $o$, the conductance $c_t(o,y)$ and the
total conductance $\pi_t(o)$ converge, with a strictly positive
limit for the latter. Thus $K_t(o,y)$ converges, and
$K_t(o,y)-K_{t+1}(o,y)\to0$. Using $0\leq u_{t+1}\leq1$ and
summing over the finitely many neighbours, we get
\[
 \sum_{y\sim o}K_t(o,y)u_{t+1}(y)
 =\sum_{y\sim o}K_{t+1}(o,y)u_{t+1}(y)+o(1)
 =\frac{\kappa_{t+1}(o)}{\pi_{t+1}(o)}+o(1).
\]
The last equality follows from the capacity--return identity for
the network frozen at time $t+1$. The capacity ratio tends to
$\kappa_*(o)/\pi_\infty(o)$ and $\eps_{t+1}(o)\to0$.
Taking the lower limit in Inequality~\eqref{eq:last-departure-average}
now gives Inequality~\eqref{eq:last-exit}.

It follows that we may choose a deterministic
$t_0=t_0(o)$ and $\eta=\eta(o)>0$ such that $q_t(o)\geq\eta$
for every $t\geq t_0$. Let $\sigma_k$ be the time of the $k$-th
visit to $o$ at or after $\max\{s,t_0\}$. On the event
$\{\sigma_k<\infty\}$, the conditional probability of a further
visit is $1-q_{\sigma_k}(o)$, by the strong Markov property of the
space-time chain. This is at most $1-\eta$, so
\[
 \PP_{s,x}(\sigma_{k+1}<\infty)
 \leq(1-\eta)\PP_{s,x}(\sigma_k<\infty),\qquad k\geq1.
\]
Iterating gives $\PP_{s,x}(\sigma_k<\infty)\leq(1-\eta)^{k-1}$
for $k\geq1$. Summing over $k$, we find that the expected number
of visits at or after $\max\{s,t_0\}$ is at most $1/\eta$.
There are at most $(t_0-s)_+$ visits before this time. Thus the
total number of visits to $o$ has finite expectation and is almost
surely finite. Taking a countable intersection over
$o\in V$, we obtain the almost-sure assertion for every vertex
simultaneously. This proves Theorem~\ref{thm:transience}.

\section{Recurrence through finite-volume plateaux}\label{sec:plateaux}

In this section, we prove Theorem~\ref{thm:plateau} and
Corollary~\ref{cor:boxes}. We first work on a fixed finite set. A lower
bound on transition probabilities and an upper bound on returns give a
lower bound on the probability of visiting a given vertex. We then apply
this bound at the beginning of each plateau, conditionally on the past.
For the application to boxes, it remains to check that the two estimates
hold uniformly as the boxes grow.

\subsection{A hitting estimate on one plateau}

\begin{lemma}\label{lem:one-plateau}
Let $P$ be a Markov transition operator on a finite set $W$, with invariant
probability measure $\widehat\pi$. Fix $z\in W$ such that
$\widehat\pi(z)>0$, and let $b\ge0$ and $u\ge1$ be integers.
Assume that there exist constants $\alpha>0$ and $B>0$ such that
\begin{align}
 P^m(x,z)&\ge \alpha\widehat\pi(z),
 &&x\in W,\quad m\ge b,\label{eq:one-plateau-lower}\\
 \sum_{k=0}^{u}P^k(z,z)&\le
 B\bigl(1+u\widehat\pi(z)\bigr).
 \label{eq:one-plateau-green}
\end{align}
Then, for every $x\in W$,
\begin{equation}\label{eq:one-plateau-hit}
 \PP_x\bigl(X_j=z\text{ for some }j\in\{b+1,\ldots,b+u\}\bigr)
 \ge \frac{\alpha}{4B}\min\{1,u\widehat\pi(z)\}.
\end{equation}
\end{lemma}
\begin{proof}
Let
\[
 N=\sum_{j=1}^{u}\1_{\{X_{b+j}=z\}}.
\]
By Inequality~\eqref{eq:one-plateau-lower}, we have
\begin{equation}\label{eq:plateau-first-moment}
 \EE_x[N]=\sum_{j=1}^{u}P^{b+j}(x,z)
 \ge \alpha u\widehat\pi(z).
\end{equation}
To bound the second moment, let us condition at the earlier of two visit
times. By the Markov property,
\begin{align*}
 \EE_x[N^2]
 &=\EE_x[N]+2\sum_{1\le i<j\le u}
     P^{b+i}(x,z)P^{j-i}(z,z)\\
 &\le\EE_x[N]\left(1+2\sum_{k=1}^{u}P^k(z,z)\right)\\
 &\le 2B\EE_x[N]\bigl(1+u\widehat\pi(z)\bigr).
\end{align*}
In the last inequality we used~\eqref{eq:one-plateau-green}, including
the term $P^0(z,z)=1$. Since $N=N\1_{\{N>0\}}$, Cauchy--Schwarz and
\eqref{eq:plateau-first-moment} give
\[
 \PP_x(N>0)\ge\frac{\EE_x[N]^2}{\EE_x[N^2]}
 \ge\frac{\alpha}{2B}\frac{u\widehat\pi(z)}{1+u\widehat\pi(z)}.
\]
Since $r/(1+r)\ge\frac12\min\{1,r\}$ for $r\ge0$, this proves
Inequality~\eqref{eq:one-plateau-hit}.
\end{proof}

\subsection{Proof of the plateau criterion}

We use the notation of Theorem~\ref{thm:plateau}. We consider finite sets
$V_1\subseteq V_2\subseteq\cdots$ with union $V_\infty$,
$T_1=0$, and, for every plateau, $L_n=T_{n+1}-T_n$ is a positive integer.
The transition from time $t$ uses $P_n$ whenever $T_n\le t<T_{n+1}$.
Thus, conditional on the history up to $T_n$, the sequence
\[
 (X_{T_n},X_{T_n+1},\ldots,X_{T_{n+1}})
\]
is a Markov chain run for $L_n$ steps with transition operator $P_n$
from its entrance point. Its last state belongs to $V_n$, and therefore
also to $V_{n+1}$. We set $\mathcal F_t=\sigma(X_0,\ldots,X_t)$ and
$u_n=(L_n-b_n)_+$.

\begin{proof}[Proof of Theorem~\ref{thm:plateau}]
Fix $z\in V_\infty$ and choose $n_0$ such that $z\in V_n$ for all $n\ge n_0$.
For $n\ge n_0$ with $u_n\ge1$, let
\[
 A_n(z)=\{X_t=z\text{ for some }
 t\in\{T_n+b_n+1,\ldots,T_n+L_n\}\}.
\]
Set $A_n(z)=\varnothing$ if $u_n=0$.
By Lemma~\ref{lem:one-plateau}, applied conditionally on
$\mathcal F_{T_n}$, we have
\begin{equation}\label{eq:plateau-conditional-success}
 \PP(A_n(z)\mid\mathcal F_{T_n})
 \ge c\min\{1,u_n\widehat\pi_n(z)\}=:q_n,
 \qquad c=\frac{\alpha}{4B}.
\end{equation}
Here $q_n$ is the same for every entrance point $X_{T_n}\in V_n$.

For $M\ge N\ge n_0$, the event
$\bigcap_{n=N}^{M-1}A_n(z)^c$ belongs to $\mathcal F_{T_M}$.
Conditioning on this sigma-field and using
Inequality~\eqref{eq:plateau-conditional-success} gives
\begin{align*}
 \PP\left(\bigcap_{n=N}^{M}A_n(z)^c\right)
 &=\EE\left[
   \1_{\bigcap_{n=N}^{M-1}A_n(z)^c}
   \PP(A_M(z)^c\mid\mathcal F_{T_M})\right]\\
 &\le(1-q_M)
   \PP\left(\bigcap_{n=N}^{M-1}A_n(z)^c\right)\\
 &\le\prod_{n=N}^{M}(1-q_n)
 \le\exp\left\{-\sum_{n=N}^{M}q_n\right\}.
\end{align*}
By assumption, the series $\sum_nu_n\widehat\pi_n(z)$ diverges.
It follows that $\sum_nq_n=\infty$. Indeed, either infinitely many terms
$u_n\widehat\pi_n(z)$ exceed one, or the two series have proportional
tails. Letting $M\to\infty$, we see that, for every $N$, the probability
of missing $z$ in all these intervals from plateau $N$ onwards is zero. Therefore,
$A_n(z)$ occurs infinitely often almost surely. The visits occur in
the disjoint integer intervals $(T_n,T_{n+1}]$, so their times are
distinct. Since $V_\infty$ is countable, we may intersect these
probability-one events over $z\in V_\infty$. This proves the theorem.
\end{proof}

The argument up to this point applies to nonreversible kernels as well.
For the reflected walk on a box we shall use reversibility to check
the two hypotheses of Lemma~\ref{lem:one-plateau}.

\subsection{Uniform estimates for reflected boxes}

For every $R\ge1$, let
\[
 B_R=\{x\in\Z^d:\|x\|_\infty\le R\}.
\]
We equip $B_R$ with its induced nearest-neighbour graph. Assume that the
symmetric edge conductances satisfy
\[
 \kappa\le c(e)\le\kappa^{-1}
\]
for some fixed $\kappa\in(0,1]$. For every vertex of this graph, define
\[
 w(x)=\sum_{y\in B_R:y\sim x}c(x,y),\qquad
 Z_R=\sum_{x\in B_R}w(x),\qquad
 \widehat\pi(x)=\frac{w(x)}{Z_R}.
\]
Let $\widetilde P(x,y)=c(x,y)/w(x)$ for neighbouring vertices of $B_R$
and let $\widetilde P(x,y)=0$ otherwise. At a boundary vertex, $w(x)$
counts only edges inside the box. We use the lazy walk
$P=(I+\widetilde P)/2$.

\begin{lemma}\label{lem:uniform-box-estimates}
There exist constants $a_0,c_0,C_0>0$, depending only on $d$ and $\kappa$,
such that every transition operator defined above satisfies
\begin{align}
 c_0R^{-d}&\le\widehat\pi(z)\le C_0R^{-d},
 &&z\in B_R,\label{eq:box-stationary-mass}\\
 P^k(z,z)&\le C_0\bigl((k+1)^{-d/2}+R^{-d}\bigr),
 &&k\ge0,\label{eq:box-diagonal}\\
 P^m(x,z)&\ge\tfrac12\widehat\pi(z),
 &&x,z\in B_R,\quad m\ge\lceil a_0R^2\rceil.
 \label{eq:box-relative-mixing}
\end{align}
\end{lemma}
\begin{proof}
Throughout the proof, the constants depend only on $d$ and $\kappa$.
First, even at the corners of the box the degree is at least $d$, and it is at most
$2d$ throughout the box. Thus $d\kappa\le w(x)\le 2d/\kappa$ and $Z_R\asymp R^d$, giving
Inequality~\eqref{eq:box-stationary-mass}.

For the other two estimates, we use Nash and Poincar\'e inequalities.
We denote the norms with respect to the unnormalized measure by
$\|\cdot\|_{p,w}$ and define
\[
 \mathcal E_P(f,f)=\langle f,(I-P)f\rangle_w.
\]
For every function on the box, the following Nash and Poincar\'e
inequalities hold:
\begin{align}
 \|f\|_{2,w}^{2+4/d}
 &\le C\bigl(\mathcal E_P(f,f)+R^{-2}\|f\|_{2,w}^2\bigr)
       \|f\|_{1,w}^{4/d},\label{eq:box-finite-nash}\\
 \|f-\widehat\pi(f)\|_{2,w}^2
 &\le CR^2\mathcal E_P(f,f).
 \label{eq:box-poincare}
\end{align}
Let us explain the reflection argument for counting measure and unit
conductances. Translate the box to $\{0,\ldots,m-1\}^d$, where
$m=2\lfloor R\rfloor+1$. In each coordinate, extend $f$ by reflection
in the order $0,1,\ldots,m-1,m-1,\ldots,1,0$ and then periodically.
Write $F$ for the extension. In one period each vertex and each edge
of the induced box graph has $2^d$ preimages. The edges across a
reflection plane contribute zero to the energy. The Poincar\'e
inequality on the torus of side $2m$, whose spectral gap is of order
$m^{-2}$, therefore gives the second bound for counting measure.

For the Nash inequality, multiply $F$ by a cutoff which is one on the
box, vanishes at distance $m$ from it and has increments bounded by
$1/m$ on nearest-neighbour edges. The support meets only a bounded
number of copies of the box, depending on $d$. Writing $D(f)$ for
the sum of squared differences over its unoriented edges, the product
has energy at most $C_d(D(f)+m^{-2}\|f\|_2^2)$ and $\ell^1$ norm at
most $C_d\|f\|_1$. Its $\ell^2$ norm is at least $\|f\|_2$, so the
usual Nash inequality on $\Z^d$ gives the first bound. Finally,
ellipticity compares $w$ to counting measure and $\mathcal E_P$ to
$D$. For the variance we use
$\|f-\widehat\pi(f)\|_{2,w}^2=\min_{a\in\R}\|f-a\|_{2,w}^2$.
This proves~\eqref{eq:box-finite-nash} and~\eqref{eq:box-poincare}
with constants depending only on $d$ and $\kappa$. The same comparison
with unit conductances is used on a torus in~\cite{Moumeni26} to
estimate merging times. For heat-kernel estimates on graphs under
volume doubling and a Poincar\'e inequality, see~\cite{Delmotte99}.

We first prove Inequality~\eqref{eq:box-diagonal}. Fix $z$ and, for each
integer $j\ge0$, let
\[
 f_j=P^j\left(\frac{\1_{\{z\}}}{w(z)}\right),\qquad
 A_j=\|f_j\|_{2,w}^2=\frac{P^{2j}(z,z)}{w(z)}.
\]
Since the functions are nonnegative and the measure is invariant, we have
$\|f_j\|_{1,w}=1$. Moreover, $P$ is self-adjoint and its spectrum lies in
$[0,1]$, so
\[
 A_j-A_{j+1}
 =\langle f_j,(I-P^2)f_j\rangle_w
 \ge\mathcal E_P(f_j,f_j).
\]
Combining \eqref{eq:box-finite-nash} with the preceding inequality gives
\[
 A_j^{1+2/d}
 \le C\bigl(A_j-A_{j+1}+R^{-2}A_j\bigr).
\]
If $A_j\ge C_1R^{-d}$, where $C_1$ is sufficiently large, the last term
on the right is at most half of the left-hand side. Consequently,
\[
 A_j-A_{j+1}\ge c_1A_j^{1+2/d}.
\]
Applying the mean-value inequality to $r\mapsto r^{-2/d}$ gives
$A_{j+1}^{-2/d}-A_j^{-2/d}\ge2c_1/d$ as long as $A_j$ stays above this
threshold. Once it falls below the threshold, it stays below it by
contraction. Since $A_0=1/w(z)$ is uniformly bounded, we obtain
\[
 A_j\le C\bigl((j+1)^{-d/2}+R^{-d}\bigr).
\]
Since the spectrum is nonnegative, $P^k(z,z)$ is nonincreasing in $k$.
Thus the last bound proves Inequality~\eqref{eq:box-diagonal} for both
even and odd times.

It remains to prove Inequality~\eqref{eq:box-relative-mixing}.
By Inequality~\eqref{eq:box-poincare}, the spectral gap is at least
$cR^{-2}$. Let $r=\lceil R^2\rceil$ and
$g_x=\1_{\{x\}}/\widehat\pi(x)-1$.
Using Inequalities~\eqref{eq:box-diagonal} and
\eqref{eq:box-stationary-mass}, we get
\[
 \|P^rg_x\|_{2,\widehat\pi}^2
 =\frac{P^{2r}(x,x)}{\widehat\pi(x)}-1\le C.
\]
Both $P^r g_x$ and $P^r g_z$ have mean zero. For $m\ge2r$, we use
self-adjointness in the first line and the spectral gap in the second:
\begin{align*}
 \left|\frac{P^m(x,z)}{\widehat\pi(z)}-1\right|
 &=\bigl|\langle P^rg_x,P^{m-2r}P^rg_z\rangle_{\widehat\pi}\bigr|\\
 &\le C\exp\{-c(m-2r)/R^2\}.
\end{align*}
We can now choose $a_0$ sufficiently large so that the right-hand
side is at most $1/2$ whenever $m\ge\lceil a_0R^2\rceil$. This proves
Inequality~\eqref{eq:box-relative-mixing}.
\end{proof}

\begin{proof}[Proof of Corollary~\ref{cor:boxes}]
By Lemma~\ref{lem:uniform-box-estimates}, the same constants apply to
every plateau operator. Let
\[
 b_n=\lceil a_0R_n^2\rceil,\qquad u_n=(L_n-b_n)_+.
\]
Inequality~\eqref{eq:box-relative-mixing} gives the lower bound required in
Theorem~\ref{thm:plateau}, with $\alpha=1/2$.
Since $d\ge3$, the series $\sum_{k\ge0}(k+1)^{-d/2}$ converges.
Summing Inequality~\eqref{eq:box-diagonal} and using
Inequality~\eqref{eq:box-stationary-mass}, we obtain,
for all $u\ge1$ and $z\in B_{R_n}$,
\begin{align*}
 \sum_{k=0}^{u}P_n^k(z,z)
 &\le C\left(1+\frac{u}{R_n^d}\right)\\
 &\le B\bigl(1+u\widehat\pi_n(z)\bigr),
\end{align*}
where $B$ does not depend on $n,z,u$.
For every fixed $z$, we have $z\in B_{R_n}$ for all sufficiently large
plateaux. By Inequality~\eqref{eq:box-stationary-mass},
\[
 \sum_n\frac{u_n}{R_n^d}=\infty
 \quad\Longrightarrow\quad
 \sum_{n:z\in B_{R_n}}u_n\widehat\pi_n(z)=\infty.
\]
Theorem~\ref{thm:plateau} now gives infinitely many visits to every
vertex almost surely.
\end{proof}

In particular, assume that $L_n\ge2b_n$ for all sufficiently large
plateaux. Then $L_n/2\le u_n\le L_n$, so the sufficient condition is
equivalent to
\begin{equation}\label{eq:plateau-box-series}
 \sum_n\frac{L_n}{R_n^d}=\infty.
\end{equation}
The comparison between $L_n$ and $b_n$ is needed to replace $(L_n-b_n)_+$
by $L_n$ in the divergence condition.

\subsection{Two choices of plateaux}

We first take $R_n=2^n$ and
\[
 L_n=\left\lceil\frac{2^{dn}}{n}\right\rceil,
 \qquad d\ge3.
\]
In this case, $L_n/b_n\to\infty$ and $L_n/R_n^d\sim1/n$.
The series in \eqref{eq:plateau-box-series} therefore diverges.
Corollary~\ref{cor:boxes} implies that every vertex is visited infinitely
often.

Now take $R_n=n$ and $L_n=\lceil Cn^{d-1}\rceil$.
For $d>3$, any $C>0$ gives $L_n/b_n\to\infty$.
For $d=3$, the ratio $L_n/b_n$ tends to $C/a_0$. Choosing
$C>2a_0$ gives $L_n\ge2b_n$ eventually.
In both cases, $L_n/R_n^d\sim C/n$. Once again the series in
\eqref{eq:plateau-box-series} diverges, and Corollary~\ref{cor:boxes}
gives recurrence.

In both examples the divergence comes from visits after the mixing
time $b_n$. When the series in~\eqref{eq:box-series} converges, the
estimate on these visits does not suffice to deduce transience.

\newpage
\appendix
\section{Electrical interpretation and an infinitesimal proof}\label{sec:electrical}

Let $(V,E)$ be a finite connected undirected graph and let
$\varnothing\ne A\subset V$. We write
\[
\mathcal F=\{f\in\R^V:f|_A=0\}.
\]
Let $g:A\to\R$ give the prescribed boundary voltages. We keep $g$
fixed while the conductances vary. For every edge $e=\{x,y\}\in E$,
choose an orientation from $x$ to $y$. For $f\in\R^V$, set $\nabla_e f=f(y)-f(x)$.
Assume that, for every $e\in E$ and every $\theta\in[0,1]$,
\[
c_\theta(e)=c_0(e)+\theta\Delta c(e),\qquad
c_0(e)>0,\qquad \Delta c(e)\ge0.
\]
For $f,h\in\R^V$, define
\[
\mathcal E_\theta(f,h)=\sum_{e\in E}c_\theta(e)
\nabla_e f\,\nabla_e h,
\qquad
\Delta\mathcal E(f,h)=\sum_{e\in E}\Delta c(e)
\nabla_e f\,\nabla_e h.
\]
The two forms do not depend on the chosen orientations. Since the graph
is connected, $A\ne\varnothing$ and all conductances are positive,
$\mathcal E_\theta$ is positive definite on $\mathcal F$.

For every $\theta\in[0,1]$, let $u_\theta$ be the unique function satisfying
\[
u_\theta|_A=g,\qquad
\forall f\in\mathcal F,\qquad
\mathcal E_\theta(u_\theta,f)=0.
\]
Equivalently, $u_\theta$ minimizes $\mathcal E_\theta(u,u)$ among all functions
$u\in\R^V$ such that $u|_A=g$.
We also fix a linear functional $\nu$ on $\mathcal F$. Its
coefficients $\nu(x)$, $x\notin A$, are the prescribed interior
currents. We use the same notation $\nu$ for the extension to
$\R^V$ given by $\nu(f)=\sum_{x\notin A}\nu(x)f(x)$. Let $\phi_\theta\in\mathcal F$ be the unique
function satisfying
\[
\forall f\in\mathcal F,\qquad
\mathcal E_\theta(\phi_\theta,f)=\nu(f).
\]
For $\theta\in[0,1]$, set
\[
C_\theta=\mathcal E_\theta(u_\theta,u_\theta),\qquad
H_\theta=\nu(\phi_\theta),\qquad M_\theta=\nu(u_\theta).
\]
For every $f\in\mathcal F$, the defining equality for $\phi_\theta$ gives
\[
2\nu(f)-\mathcal E_\theta(f,f)
=\mathcal E_\theta(\phi_\theta,\phi_\theta)
-\mathcal E_\theta(f-\phi_\theta,f-\phi_\theta).
\]
Since $\mathcal E_\theta(f-\phi_\theta,f-\phi_\theta)\geq0$,
the maximum is attained at $f=\phi_\theta$. It follows that
\begin{equation}\label{eq:electrical-dual-energy}
H_\theta=\mathcal E_\theta(\phi_\theta,\phi_\theta)
=\sup_{f\in\mathcal F}
\{2\nu(f)-\mathcal E_\theta(f,f)\}.
\end{equation}

\begin{proposition}\label{prop:electrical-interpolation}
The functions $u_\theta$, $\phi_\theta$, $C_\theta$, $H_\theta$ and
$M_\theta$ are continuously differentiable in $\theta$. Moreover, for
every $\theta\in[0,1]$, the following equalities hold:
\begin{align}
C_\theta'&=\Delta\mathcal E(u_\theta,u_\theta),
\label{eq:electrical-C-derivative}\\
H_\theta'&=-\Delta\mathcal E(\phi_\theta,\phi_\theta),
\label{eq:electrical-H-derivative}\\
M_\theta'&=-\Delta\mathcal E(u_\theta,\phi_\theta).
\label{eq:electrical-M-derivative}
\end{align}
In particular,
\begin{equation}\label{eq:electrical-two-form-bound}
|M_1-M_0|\le
\sqrt{(C_1-C_0)(H_0-H_1)}.
\end{equation}
\end{proposition}
\begin{proof}
Since $V$ is finite, the equations defining $u_\theta$ and $\phi_\theta$
are finite linear systems with the same Dirichlet matrix. This
matrix is invertible on $\mathcal F$ and depends affinely on
$\theta$. Its inverse is continuously differentiable on a neighbourhood
of $[0,1]$. The right-hand sides are affine in $\theta$, so both
potentials, and hence $C_\theta$, $H_\theta$ and $M_\theta$, are
continuously differentiable up to the endpoints.

Fix $\theta\in[0,1]$. Since the boundary values do not change with
$\theta$, we have $u_\theta',\phi_\theta'\in\mathcal F$.
Differentiating the defining equations gives, for every $f\in\mathcal F$,
\begin{align*}
\mathcal E_\theta(u_\theta',f)
&=-\Delta\mathcal E(u_\theta,f),\\
\mathcal E_\theta(\phi_\theta',f)
&=-\Delta\mathcal E(\phi_\theta,f).
\end{align*}
Since $u_\theta$ is harmonic and $u_\theta'\in\mathcal F$, we have
$\mathcal E_\theta(u_\theta,u_\theta')=0$. Therefore,
\[
C_\theta'
=\Delta\mathcal E(u_\theta,u_\theta)
+2\mathcal E_\theta(u_\theta,u_\theta')
=\Delta\mathcal E(u_\theta,u_\theta).
\]
For the derivative of $H_\theta$, take $f=\phi_\theta$ in the
second differentiated equation. This gives
\[
H_\theta'=\nu(\phi_\theta')
=\mathcal E_\theta(\phi_\theta,\phi_\theta')
=-\Delta\mathcal E(\phi_\theta,\phi_\theta).
\]
Taking $f=\phi_\theta$ in the first differentiated equation instead gives
\[
M_\theta'=\nu(u_\theta')
=\mathcal E_\theta(\phi_\theta,u_\theta')
=-\Delta\mathcal E(u_\theta,\phi_\theta).
\]

It remains to prove Inequality~\eqref{eq:electrical-two-form-bound}.
Since $\Delta c(e)\ge0$ for every $e\in E$, the form
$\Delta\mathcal E$ is positive semidefinite. By Cauchy--Schwarz and
Equalities~\eqref{eq:electrical-C-derivative}--\eqref{eq:electrical-M-derivative},
\[
|M_\theta'|
\le\sqrt{\Delta\mathcal E(u_\theta,u_\theta)
\Delta\mathcal E(\phi_\theta,\phi_\theta)}
=\sqrt{C_\theta'(-H_\theta')}.
\]
Integrating and applying Cauchy--Schwarz on $[0,1]$, we get
\begin{align*}
|M_1-M_0|
&\le\int_0^1
\sqrt{\Delta\mathcal E(u_\theta,u_\theta)
\Delta\mathcal E(\phi_\theta,\phi_\theta)}\,d\theta\\
&\le\left(\int_0^1 C_\theta'\,d\theta\right)^{1/2}
\left(\int_0^1(-H_\theta')\,d\theta\right)^{1/2}.
\end{align*}
The two integrals equal $C_1-C_0$ and $H_0-H_1$, respectively.
This proves Inequality~\eqref{eq:electrical-two-form-bound}.
\end{proof}

Let us interpret these identities in terms of the electrical network.
The quantity $C_\theta$ is the power dissipated when the boundary voltage is fixed at $g$.
Equality~\eqref{eq:electrical-C-derivative} shows that it is nondecreasing.
The potential $\phi_\theta$ corresponds to the fixed interior currents
$\nu$, with voltage zero on $A$. Its dissipated power is $H_\theta$,
which is nonincreasing by Equality~\eqref{eq:electrical-H-derivative}.

To see where the minus sign in $H_\theta'$ comes from, recall that
with fixed currents the functional maximized at $\phi_\theta$ is $f\mapsto2\nu(f)-\mathcal E_\theta(f,f)$, as in
\eqref{eq:electrical-dual-energy}. To differentiate the energy, one must
also account for the change in potential. Indeed,
\begin{align*}
\frac{d}{d\theta}\mathcal E_\theta(\phi_\theta,\phi_\theta)
&=\Delta\mathcal E(\phi_\theta,\phi_\theta)
+2\mathcal E_\theta(\phi_\theta,\phi_\theta'),\\
2\mathcal E_\theta(\phi_\theta,\phi_\theta')
&=-2\Delta\mathcal E(\phi_\theta,\phi_\theta).
\end{align*}
Finally, $M_\theta=\nu(u_\theta)$ is the pairing of the interior currents
with the voltage $u_\theta$. It should not be identified with
$\mathcal E_\theta(u_\theta,\phi_\theta)$, which is zero because
$\phi_\theta\in\mathcal F$ and $u_\theta$ is harmonic.
In Equality~\eqref{eq:electrical-M-derivative}, each increment
$\Delta c(e)$ is multiplied by the product of the two voltage drops
across $e$. This product may have either sign. Inequality~\eqref{eq:electrical-two-form-bound}
therefore controls the size of the change in $M_\theta$ through the
increase of $C_\theta$ and the decrease of $H_\theta$.

\newpage
\bibliographystyle{alpha}
\bibliography{references}
\end{document}